\documentclass[11 pt]{amsart}
\usepackage{amsmath,amssymb,amsthm}
\newtheorem{Theorem}{Theorem}[section]

\newtheorem{Corollary}[Theorem]{Corollary}

\theoremstyle{definition}
\newtheorem{Example}[Theorem]{Example}
\newtheorem{Remark}[Theorem]{Remark}

\begin{document}   


\title[cohomology of invariant subspaces]
{On cohomology of invariant submanifolds of Hamiltonian actions}

\author{Y{\i}ld{\i}ray Ozan}
\address{Department of Mathematics, Middle East Technical University,
06531 \newline Ankara, TURKEY} \email{ozan@metu.edu.tr}
\date{\today}
\thanks{The author is partially supported by the Turkish Academy of
Sciences (TUBA-GEBIP-2004-17).}  \pagenumbering{arabic}

\maketitle

\section{Introduction}
In \cite{Oz2} the author proved that if there is a free algebraic
circle action on a nonsingular real algebraic variety $X$ then the
fundamental class is trivial in any nonsingular projective
complexification  $i:X\rightarrow X_{\mathbb C}$.   The K\"ahler
forms on ${\mathbb C}^N$ and ${\mathbb CP}^N$ naturally induce
symplectic structures on complex algebraic affine or projective
varieties and in case they are defined over reals, their real
parts, if not empty, are Lagrangian submanifolds.

The following result can be considered as a symplectic equivalent
of author's above result on real algebraic varieties.

\begin{Theorem}\label{thm-B}
Assume that $G$ is $S^1$ or $S^3$ acting on a compact symplectic
manifold $(M,\omega)$ in a Hamiltonian fashion and $L^l$ is an
invariant closed submanifold. If the $G$-action on $L$ is locally
free then the homomorphism induced by the inclusion,
$i:L\rightarrow M$,
$$H_i(L,{\mathbb Q})\rightarrow H_i(M,{\mathbb Q})$$ is trivial
for $i\geq l-k+1$, where $k=\dim (G)$. In particular, the
fundamental class $[L]$ is trivial in $H_l(M,{\mathbb Q})$.

Moreover, if the corresponding sphere bundle $S^k\rightarrow
L\times EG\rightarrow L_G$  has non torsion Euler class then the
homomorphism $$i_*:H_{l-k}(L,{\mathbb Q})\rightarrow
H_{l-k}(M,{\mathbb Q}),$$ induced by the inclusion $i:L\rightarrow
M$, is also trivial (see Section 2 for the definition of $EG$ and
$L_G$).
\end{Theorem}

Since any compact connected Lie group has a circle subgroup we
deduce the following immediate corollary.

\begin{Corollary} \label{cor-A}
Let $G$ be a compact connected Lie group acting on a compact
symplectic manifold $(M,\omega)$ in a Hamiltonian fashion and $L$
is an invariant closed submanifold of dimension $l$. If the
$G$-action on $L$ is locally free then the fundamental class $[L]$
is trivial in $H_l(M,{\mathbb Q})$.
\end{Corollary}

\begin{Remark}\label{rem-A}
{\bf 1)} It is well known that the natural actions of $U(n)$ and
$T^{n}$ on the complex projective space ${\mathbb CP}^{n-1}$ and hence
on any smooth projective variety, regarded as a symplectic
manifolds with their Fubini-Study forms, are Hamiltonian (cf. see
p.163 in {\cite{McS}}).

{\bf 2)} Consider the $2$-torus $T^2=S^1\times S^1$ with the
symplectic (volume) form $d\theta_1 \wedge d\theta_2$.  Then, the
$S^1$ action on \ $T^2$, given by \ $z\cdot (w_1,w_2)=(z\cdot
w_1,w_2)$ \ is clearly symplectic but not Hamiltonian (it has no
fixed point). Let $L$ be the invariant submanifold
$S^1\times\{pt\}$, on which the circle action is free.  Clearly,
the homology class  \ $[L]$ \ is not zero in \ $H_1(T^2,{\mathbb
Q})$.  Hence, the assumption in the above theorem that the action
is Hamiltonian, is necessary.

{\bf 3)} Since $S^3=SU(2)$ is semisimple any symplectic
$SU(2)$-action is Hamiltonian (cf. see p.159 of {\cite{McS}}).
\end{Remark}

\begin{Example}
Let $G$ be a compact Lie group acting linearly on a closed
manifold $M$.  Dovermann and Masuda proved that if the action is
semifree or $G=S^1$ then there exists a nonsingular real algebraic
variety $X$ with an algebraic $G$-action equivariantly
diffeomorphic to $M$ (cf see \cite{Dov1}). If the linear action on
$X$ extends to some nonsingular projective complexification
$X_{\mathbb C}$, then by part (1) of the above remark the action
will be Hamiltonian and thus the above results can be applied to
the pair $X\subseteq X_{\mathbb C}$.
\end{Example}

Symplectic reduction and the proofs of the above results give a
somewhat stronger result.

Let $G=G_1\times \cdots \times G_d$, where each $G_j$ is either
$S^1$ or $S^3=SU(2)$ and suppose that it acts in a Hamiltonian
fashion on a closed symplectic manifold $M$, with moment map $\mu
:M \rightarrow \mathfrak{g}^*$, where $$\mu =(\mu_1, \ldots ,
\mu_d):M \rightarrow (\mathfrak{g}_1^*,\ldots,
\mathfrak{g}_d^*),$$ each $\mathfrak{g}_j^*$ being the dual of the
Lie algebra of $G_j$. Assume that $L$ is an invariant submanifold
contained in a level set $M^0= \mu^{-1}(v_1,\cdots,v_d)$ of the
moment map. Further assume that, we can form successive symplectic
reductions first by $G_1$, for level set $\mu^{-1}(v_1)$, then by
$G_2$, for the level set $(\mu_1^{-1}(v_1)\cap
\mu_2^{-1}(v_2))/G_1$, and so on for all $G_j$ (note that we use
the same notation for the moment map on the reduced spaces). These
assumptions clearly imply that the $G$-action on $L$ is locally
free.

\begin{Theorem}\label{thm-C}
Assume the above setup. Then the induced homomorphism
$i^*:H_i(L,{\mathbb Q})\rightarrow H_i(M,{\mathbb Q})$ is trivial,
for $i\geq l-k+1$, where $k=\dim (G)$ and $i:L\rightarrow M$ is
the inclusion map.
\end{Theorem}

\section{Proofs}
A $G$-space $X$ is called equivariantly formal if its equivariant
cohomology is isomorphic to its usual cohomology tensored by the
cohomology of the classifying space.  Equivariant formality
implies that the equivariant cohomology of $X$ injects into the
equivariant cohomology of the fixed point set,
$H^*_G(X)\rightarrowtail H^*_G(X^G)$ (cf. see Theorem 11.4.5 in
\cite{GS}).  On the other hand, Kirwan proved in \cite{Kir1} that
a compact Hamiltonian space $M$ is equivariantly formal.
Therefore, Theorem~\ref{thm-B} is a consequence of
Theorem~\ref{thm-A} below and the Kirwan's result.

Let $G$ act freely on $EG=S^{\infty}$, which is $\lim S^{2n-1}$ in
case $G=S^1$ and $\lim S^{4n-1}$ in case $G=SU(2)$.   Also let
$BG=EG/G$, which is ${\mathbb CP}^{\infty}$ if $G=S^1$ and
${\mathbb HP}^{\infty}$ if $G=SU(2)$. For any $G$-space $X$, we
will denote the twisted product $X\times_G EG$ by $X_G$, where
$G$-action on $X\times EG$ is given by $g\cdot (x,h)=(g^{-1}\cdot
x, h\cdot g)$, for any $g\in G$, $h \in EG$ and $x\in X$.  Then
for any coefficient ring $\it R$, the $G$-equivariant cohomology
of $X$ is defined to be the ordinary cohomology of $X_G$:
$$H_G^*(X,{\it R})\doteq H^*(X_G,{\it R}).$$ Let $p:X\times
EG\rightarrow X_G$ be the quotient map. Since $S^{\infty}$ is
contractible the induced homomorphism by $p$ in cohomology can be
regarded as a map $p^*:H_G^*(X,{\it R})\rightarrow H^*(X,{\it
R})$.

\begin{Theorem} \label{thm-A}
Let $M$ be an orientable closed manifold equipped with an
equivariantly formal $G$-action, and $L^l$ an invariant closed
submanifold. If the $G$-action on $L$ is locally free then the
homomorphism induced by the inclusion, $i:L\rightarrow M$,
$$H_i(L,{\mathbb Q})\rightarrow H_i(M,{\mathbb Q})$$ is trivial
for $i\geq l-k+1$, where $k=\dim (G)$. In particular, the
fundamental class $[L]$ is trivial in $H_l(M,{\mathbb Q})$.

Moreover, if the corresponding sphere bundle $S^k\rightarrow
L\times E_G\rightarrow L_G$  has non torsion Euler class then the
homomorphism $$i_*:H_{l-k}(L,{\mathbb Q})\rightarrow
H_{l-k}(M,{\mathbb Q}),$$ induced by the inclusion $i:L\rightarrow
M$, is also trivial.
\end{Theorem}

The above theorem is a consequence of a more general result which
can be stated using equivariant cohomology.

\begin{Theorem} \label{thm-main}
Let $M$ be an orientable closed manifold with an equivariantly
formal $G$-action and $L^l$ an invariant closed submanifold. If
$i:L\rightarrow M$ is the inclusion map then the image of
$i^*:H^*(M,{\mathbb Q})\rightarrow H^*(L,{\mathbb Q})$ lies in the
image of $p^*:H_G^*(L,{\mathbb Q})\rightarrow H^*(L,{\mathbb Q})$.
\end{Theorem}

If the $G$-action on $X$ is free then $X_G\rightarrow X/G$ has
contractible fibers and thus $X_G\rightarrow X/G$ is a homotopy
equivalence. If the action is locally free then the fibers are
finite cyclic quotients of contractible spaces and therefore the
rational cohomology of $X_G$ and $X/G$ are still isomorphic.
Hence, we obtain the following corollary.

\begin{Corollary} \label{cor-main}
Assume that $M$ and $G$ are as in the above theorem and the
$G$-action on $L$ is locally free. If $p:L\rightarrow B=L/G$ is
the quotient map, then the image of $i^*:H^*(M,{\mathbb
Q})\rightarrow H^*(L,{\mathbb Q})$ lies in the image of
$p^*:H^*(B,{\mathbb Q})\rightarrow H^*(L,{\mathbb Q})$.
\end{Corollary}

\begin{proof}[Proof of Theorem~\ref{thm-main}]
Consider the following commutative ladder of Gysin sequences
corresponding to the sphere bundles $S^k=G\rightarrow M\times
EG\rightarrow M_G$ and $S^k=G\rightarrow L\times EG\rightarrow
L_G$ \ ($k=\dim(G)$):

\begin{center}
$\cdots\rightarrow H^i(M_G, {\mathbb
Q})\stackrel{p^*}{\rightarrow} H^i(M, {\mathbb Q})
\stackrel{p^!}{\rightarrow} H^{i-k}(M_G, {\mathbb
Q})\stackrel{\cup e}{\rightarrow} H^{i+1}(M_G,{\mathbb
Q})\rightarrow \cdots $ \vspace{0.2cm}

\hspace{0.3cm}$\downarrow i^*$\hspace{1.8cm}$\downarrow i^*
$\hspace{2cm}$\downarrow i^*$\hspace{2cm}$\downarrow i^*$

\vspace{0.2cm} $\cdots\rightarrow H^i(L_G, {\mathbb Q})
\stackrel{p^*}{\rightarrow} H^i(L, {\mathbb
Q})\stackrel{p^!}{\rightarrow} H^{i-k}(L_G, {\mathbb
Q})\stackrel{\cup e}{\rightarrow} H^{i+1}(L_G,{\mathbb
Q})\rightarrow \cdots$,
\end{center}
where \ $p^!$ \ is the connecting homomorphism (can be thought as
integration along fiber) and \ $e\in H^{k+1}(M_G, {\mathbb Q})$ \
is the image of the Euler class of the sphere bundle of under the
natural map $H^{k+1}(M_G, {\mathbb Z})\rightarrow H^{k+1}(M_G,
{\mathbb Q})$.

Note that to prove the theorem it suffices to show that the map
$p^!$ in the top row is trivial.  Indeed, we claim that the map \
$H^{i-k}(M_G, {\mathbb Q})\stackrel{\cup  e}{\rightarrow}
H^{i+1}(M_G,{\mathbb Q})$ \ is injective.  To see this let \ $M^G$
\ denote the fixed point and consider the following commutative
diagram:

\begin{center}
$H_G^{i-k}(M, {\mathbb Q})\stackrel{\cup  e}{\rightarrow}
H_G^{i+1}(M,{\mathbb Q})$ \vspace{0.2cm}

$\downarrow i^*$\hspace{2cm}$\downarrow i^*$

\vspace{0.2cm} $H_G^{i-k}(M^G, {\mathbb Q})\stackrel{\cup
e}{\rightarrow} H_G^{i+1}(M^G,{\mathbb Q})$.
\end{center}
By assumption the vertical arrows are injections and hence it is
enough to show that the bottom row is injective. For the latter,
note that the $G$-action on \ $M^G$ \ is trivial and hence the
corresponding \ $G$-bundle for \ $M^G$ \ is
$$G\rightarrow M^G\times S^{\infty}\stackrel{id \times
h}{\longrightarrow} M^G\times BG,$$ where \ $BG={\mathbb
CP}^{\infty}$ or ${\mathbb HP}^{\infty}$, depending on whether \
$G$ \ is \ $S^1$  \ or \ $SU(2)=S^3$. Moreover, the Euler class of
the bundle is \ $e=(1,e_0) \in  H^0(M^G,{\mathbb Q})\times
H^{k+1}(BG,{\mathbb Q})$, where $e_0$ is a generator of
$H^{k+1}(BG,{\mathbb Q})$ \ and hence cup product with the Euler
class is injective.
\end{proof}

\begin{proof}[Proof of Theorem~\ref{thm-A}]
By the Universal Coefficient Theorem it suffices to show that the
map $i^*:H^i(M,{\mathbb Q})\rightarrow H^i(L,{\mathbb Q})$ is
trivial for $i\geq l-k+1$.  Therefore, by Corollary~\ref{cor-main}
it is enough to show that the map $$p^*:H^i(L_G,{\mathbb
Q})\rightarrow H^i(L,{\mathbb Q})$$ is trivial.  However, since
the $G$-action on $L$ is locally free the rational (co)homology of
$L_G$ is equal to that of the $l-k$-dimensional orbifold $B=L/G$.
In particular, $H^{i}(L_G, {\mathbb Q})=0$ for $i\geq l-k+1$. This
finishes the proof of the first statement.

For the second statement consider the Gysin sequence corresponding
to the $G$-bundle $p:L\times EG\rightarrow L_G$:
$$\rightarrow H^{l-2k-1}(L_G, {\mathbb Q})\stackrel{\cup
e}{\rightarrow} H^{l-k}(L_G, {\mathbb
Q})\stackrel{p^*}{\rightarrow} H^{l-k}(L, {\mathbb
Q})\stackrel{p^!}{\rightarrow} H^{l-k-1}(L_G,{\mathbb
Q})\rightarrow.$$ Since $e\in H^{k+1}(B,{\mathbb Z})$ is not a
torsion class, by Poincar\'e duality the map given by the cup
product with the Euler class is onto.  This implies that the map
$p^*$ is trivial and hence the proof concludes as in the first
statement.
\end{proof}

\begin{proof}[Proof of Theorem~\ref{thm-C}]
First we will prove that $$Im(H^i(M,{\mathbb Q})\rightarrow
H^i(L,{\mathbb Q})) \subseteq Im(H^i(L/G,{\mathbb Q})\rightarrow
H^i(L,{\mathbb Q})).\hspace{1cm}$$ Proof is by induction on $d$,
the number of factors in the decomposition $G=G_1\times \cdots
\times G_d$. The case $d=1$ is contained in
Theorem~\ref{thm-main}. Suppose that the theorem holds for all
integers $1,\ldots,d-1$, where $d\geq 2$. Consider the action of
$G_1$ on $M$, with the moment map $\mu_1 :M \rightarrow
\mathfrak{g}_1^*$, where $$\mu =(\mu_1, \ldots , \mu_d):M
\rightarrow (\mathfrak{g}_1^*,\ldots, \mathfrak{g}_d^*).$$

Let $M_{red}$ denote the reduced space $M^1=\mu_1^{-1}(v_1)/G_1$.
Note that $G_{red}=G_2\times \cdots \times G_d$ has an induced
Hamiltonian action on the symplectic manifold $M_{red}$ and the
moment map $\mu$ descends to a moment map $$\mu_{red} =(\mu_2,
\ldots , \mu_d):M_{red} \rightarrow (\mathfrak{g}_2^*,\ldots,
\mathfrak{g}_d^*),$$ satisfying the same hypothesis as $\mu$.
Abusing the notation further, we will denote $L/G_1$ by $L_{red}$.
Note that $L_{red}$ is an invariant submanifold of $M_{red}$.

Hence, by the induction hypothesis
$$\hspace{-2cm}Im(H^i(M_{red},{\mathbb Q})\rightarrow
H^i(L_{red},{\mathbb Q})) $$ $$ \hspace{3cm} \subseteq
Im(H^i(L_{red}/G_{red},{\mathbb Q})\rightarrow
H^i(L_{red},{\mathbb Q})).\hspace{1cm} (*)$$

Now, we will consider a ladder of exact sequences similar to the
one used in the proof of Theorem~\ref{thm-main}:

\begin{center}
$\cdots\rightarrow H^i(M_G, {\mathbb
Q})\stackrel{p^*}{\rightarrow} H^i(M, {\mathbb Q})
\stackrel{p^!}{\rightarrow} H^{i-k_1}(M_G, {\mathbb
Q})\stackrel{\cup e}{\rightarrow} H^{i+1}(M_G,{\mathbb
Q})\rightarrow $ \vspace{0.2cm}

\hspace{0.3cm}$\downarrow \kappa$\hspace{1.8cm}$\downarrow i^*
$\hspace{2cm}$\downarrow \kappa$\hspace{2cm}$\downarrow i^*$

\vspace{0.2cm} $\cdots\rightarrow H^i(M_{red}, {\mathbb Q})
\stackrel{p^*}{\rightarrow} H^i(M^1, {\mathbb
Q})\stackrel{p^!}{\rightarrow} H^{i-k_1}(M_{red}, {\mathbb
Q})\stackrel{\cup e}{\rightarrow} H^{i+1}(M_{red},{\mathbb
Q})\rightarrow$ \vspace{0.2cm}

\hspace{0.3cm}$\downarrow i^*$\hspace{1.8cm}$\downarrow i^*
$\hspace{2cm}$\downarrow i^*$\hspace{2cm}$\downarrow i^*$

\vspace{0.2cm} $\cdots\rightarrow H^i(L_{red}, {\mathbb Q})
\stackrel{p^*}{\rightarrow} H^i(L, {\mathbb
Q})\stackrel{p^!}{\rightarrow} H^{i-k_1}(L_{red}, {\mathbb
Q})\stackrel{\cup e}{\rightarrow} H^{i+1}(L_{red},{\mathbb
Q})\rightarrow$,
\end{center}
where $k_1=\dim (G_1)$ and the maps from the top row to the middle
one denoted by $\kappa$ and induced by inclusion maps also, are
the Kirwan maps (\cite{Kir1}). As in the proof of
Theorem~\ref{thm-main} the map $p^!$ in the top row is trivial.
Noting that $L_{red}/G_{red}=L/G$ it follows from $(*)$ and the
above diagram that

$$Im(H^i(M,{\mathbb Q})\rightarrow
H^i(L,{\mathbb Q})) \subseteq Im(H^i(L/G,{\mathbb Q})\rightarrow
H^i(L,{\mathbb Q})).$$

Finally, the arguments in the first paragraph of the proof of
Theorem~\ref{thm-A} finishes the proof.
\end{proof}

\begin{Remark}\label{rem-Kir}
It is known that the Kirwan map $\kappa$ is surjective
(\cite{Kir1}). Even though we don't need this information for the
above proof a diagram chase in the exact sequences implies the
following corollary.
\end{Remark}

\begin{Corollary}\label{cor-final}
Let $M$, $M_{red}$, $L$ and $L_{red}$ be as above.  Then the map
$$p^*:Im(H^i(M_{red},{\mathbb Q})\rightarrow H^i(L_{red},{\mathbb Q}))
\rightarrow Im(H^i(M,{\mathbb Q})\rightarrow H^i(L,{\mathbb Q}))$$
is onto for any $i$ and is an isomorphism for $i=1$.
\end{Corollary}

\subsection*{Acknowledgment} Some part of this research has been
completed during a visit to Universit\'{e} de Rennes 1, Rennes
France, in October 2002. I am grateful to Mathematics Department
for the invitation and the warm hospitality.  I am also grateful
to Selman Akbulut, Heiner Dovermann and Yael Karshon for their
comments on the earlier version of this note.  I would like thank
also Frances Kirwan for sending a copy of her book \cite{Kir1} to
me to learn more about her results.

\providecommand{\bysame}{\leavevmode\hboxto3em{\hrulefill}\thinspace}

\end{document}